\input amstex
\documentstyle{amsppt}
\NoBlackBoxes
\loadbold
\topmatter
\title
Contact deformations of closed 1-forms on $\text{\bf T}^2$ bundles over $\text{\bf S}^1$
\endtitle
\rightheadtext{Contact deformations}
\author Hamidou Dathe and Philippe Rukimbira
\endauthor
\leftheadtext{Dathe and Rukimbira}
\address Department of Mathematics, UCAD, Dakar, Senegal
\endaddress
\address Department of Mathematics, Florida International University, Miami Florida 33199, USA
\endaddress
\abstract \magnification\magstep1 \baselineskip 12pt If a closed
3-manifold $M$ supports a closed, nonsingular, irrational 1-form
which linearly deforms into contact forms, then $M$ supports a
K-contact form. On the 3-torus, a closed nonsingular 1-form deforms
linearly into contact forms if and only if it is a fibration 1-form.
On any other 2-dimensional torus bundle over the circle, every
closed, nonsingular 1-form deforms linearly into contact forms.
\endabstract
\endtopmatter
\document
\baselineskip 24pt
\head 0 Introduction
\endhead
The main result of this paper, Theorem 4.1, is a complete
classification of closed, nonsingular 1-forms using linear
deformations into contact forms on $\text{\bf T}^2$ bundle over
$\text{\bf S}^1$. We prove that every closed, nonsingular 1-form on
a $\text{\bf T}^2$ bundle other than $\text{\bf T}^3$, deforms
linearly into contact forms. On $\text{\bf T}^3$, a closed
nonsingular 1-form deforms linearly into contact forms if and only
if it is rational.

The proof relies on several results, the most crucial of which are
contained in Proposition 1.1, Proposition 3.1 and Theorem 3.1.

All results cited in this paper have been paraphrased so as to fit
our needs for clarity of presentation. We refer to their original
publications for complete statements.

\head 1 Preliminaries
\endhead
A contact form on a $2n+1$-dimensional manifold $M$ is a $1$-form
$\alpha$ such that $\alpha\wedge d\alpha$ is a volume form on $M$.
We call $(M, \alpha)$ a contact manifold. In the case where $M$ is
oriented, a positive contact form is a $1$-form such that
$\alpha\wedge (d\alpha)^{n} > 0$.
 If $\alpha\wedge (d\alpha)^{n}\ge 0$ holds, we say that $\alpha$ is a positive
 confoliation form.
If $\alpha$ is a contact form, then the system of equations
$\alpha(Z)= 1$ and $d\alpha(Z, X)= 0$ for arbitrary $X$, uniquely
determines a vector field $Z$ called the Reeb vector field, or the
characteristic vector field of $\alpha$. Using Seiberg-Witten
theory, Taubes (\cite{TAU}) has shown the following.

\proclaim{Proposition 1.1} Suppose $M$ is a closed 3-dimensional
orientable manifold and $\alpha$ is a contact form on $M$. Then the
Reeb vector field $Z$ of $\alpha$ has at least 1 periodic orbit.
\endproclaim

The tangent subbundle $\xi= Kern~ \alpha$ of rank $2n$ is called the
contact structure associated with $\alpha$. If $\alpha$ is a
confoliation form, $\xi=Kern ~\alpha$ is called the confoliation
structure associated with $\alpha$. In general, a contact structure
on a $2n+1$-dimensional manifold is a rank $2n$ tangent subbundle
which is locally determined by contact forms (see Blair's book
\cite{BLA} for details about contact structures). The manifolds in
this paper will be oriented and all the planes fields considered
herein are supposed to be transversely orientable.

\noindent Also, a contact manifold $(M, \alpha)$ with Reeb vector
field $Z$ admits a riemannian metric $g$, called a contact metric,
and a $(1,1)$ tensor field $J$ such that the following relations
holds (\cite{BLA}).

$JZ= 0$, $\alpha(Z)= 1$, $J^{2}= -I+ \alpha\otimes Z$, $\alpha(X)=
g(Z, X)$,

$g(JX, JY)= g(X, Y)- \alpha(X)\alpha(Y)$, $g(X, JY)= d\alpha(X, Y)$.

If the Reeb vector field $Z$ is Killing relative to a contact metric
$g$, then $\alpha$ is called a K-contact form. When all orbits of
$Z$ are periodic, the contact form is said to be quasi-regular. The
following result about quasi-regular contact forms was obtained in
\cite{RUK}.

\proclaim{Proposition 1.2} Let $(M,\alpha )$ be a closed, contact
manifold such that $\alpha$ is quasi-regular. Then $\alpha$ is a
K-contact form.\endproclaim

 Another result about K-contact forms is taken from \cite{RU1}.
 \proclaim{Proposition 1.3} No torus can carry a K-contact form.
 \endproclaim

Let $\xi$ a plane field on a manifold $M$. When $\xi$ is a
foliation, we say that $\xi$ is  deformable into contact  structures
if there exists a one parameter family $\xi_{t}$ of hyperplanes
fields satisfying $\xi_{0}= \xi$ and, for all $t>0$, $\xi_{t}$ is
contact. In this paper, we deal with particular deformations called
{\it linear}. For a foliation defined by a $1$-form $\alpha_{0}$, a
deformation into contact structures $\xi_{t}$ defined by $1$-forms
$\alpha_{t}$ is said to be linear if $\alpha_{t}= \alpha_{0}+
t\alpha$ where $\alpha$ is a $1$-form (independent of $t$).

\head 2. Linear deformations into contact forms\endhead

 \noindent In
(\cite{DA1}), we studied linear deformations of closed nonsingular
1-forms and observed the fact that not any closed, nonsingular
1-form deforms linearly into contact forms. The following result was
obtained as a characterization of those closed, nonsingular 1-forms
which admit linear deformations into contact forms.

\proclaim{Proposition 2.1}On a closed $2n+1$-dimensional manifold
with a closed, nonsingular 1-form $\alpha_0$, the family of 1-forms
$\alpha_t=\alpha_0+t\alpha$ is contact for all $t>0$ if and only if
$\alpha$ is a contact form and $\alpha_0(Z)=0$ where $Z$ is the Reeb
vector field of $\alpha$.
\endproclaim

Proposition 2.1 was later used by the authors in (\cite{DA2}) to
obtain a characterization of $\text{\bf T}^2$ bundles over
$\text{\bf S}^1$.

\proclaim{Proposition 2.2} A closed 3-dimensional manifold $M$ is a
$\text{\bf T}^2$ bundles over $\text{\bf S}^1$ if and only if it
carries a nonsingular, closed 1-form that is linearly deformable
into contact forms.\endproclaim

\head 3. Isotopy of nonsingular 1-forms\endhead \noindent Two forms
$\omega_1$ and $\omega_2$ on a manifold $M$ are said to be isotopic
if there exist a diffeomorphism $\phi$ isotopic to the identity such
that $\omega_2=\phi^\ast \omega_1$ the pullback of $\omega_1$. The
following result is due to Blank and Laudenback (\cite{LAB})

\proclaim{Proposition 3.1} Suppose $M$ is a closed 3-dimensional
manifold and $\omega_1$ and $\omega_2$ are two cohomologous, closed,
nonsingular 1-forms on $M$. Then $\omega_1$ and $\omega_2$ are
isotopic.\endproclaim

For contact forms, we prove the following isotopy result.

\proclaim{Proposition 3.2} Suppose, on a closed, $2n+1$-dimensional
manifold with a closed, nonsingular 1-form $\alpha_0$, that the
family of 1-forms $\alpha_t=\alpha_0+t\alpha$ is contact for all
$t>0$. Then for each $t>0$, $t\alpha$ and $\alpha_t$ are isotopic
through diffeomorphisms preserving the Reeb foliation of $\alpha$.
\endproclaim

\demo{Proof} Let $\beta_s=\alpha_t-s\alpha_0$. A direct calculation
shows that $\beta_s$ is contact for for all $s$. Notice that
$\beta_s$ and $\alpha_t$ have the same Reeb field ${Z\over t}$.
Moreover, $\beta_0=\alpha_t$, $\beta_1=t\alpha$ and
${d\over{ds}}\beta_s=\alpha_0$. Since (by Proposition 2.1)
$\alpha_0({Z\over t})=0$ and $i_{Z\over t}d\beta_s=0$, there is a
unique infinitesimal isotopy $X_s$ determined by the two equations:
$\beta_s(X_s)=0$ and $i_{X_s}d\beta_s=-\alpha_0$. Denote by
$\Theta_s$ the isotopy generated by $X_s$. Since
$\Theta_0^\ast\beta_0=\alpha_t$ and
${d\over{ds}}(\Theta^\ast_s\beta_s)=\Theta^\ast_s(L_{X_s}\beta_s+\alpha_0)=
\Theta^\ast_s(i_{X_s}d\beta_s+\alpha_0)=0$, it follows that
$$\Theta^\ast_s\beta_s=\alpha_t$$ for all $s$. So one obtains the
identity $\alpha_t=\Theta^\ast_1\beta_1=\Theta^\ast_1(t\alpha
)$.That each $\Theta_s$ preserves the Reeb foliation of $\alpha$
follows from the fact that $\beta_s$ and $t\alpha$ share the same
Reeb vector field ${Z\over t}$. \quad\qed
\enddemo

A contact structure $kern~\alpha$ on a 3-manifold $M$ is said to be
overtwisted (see \cite{ELI}) if there exists an embedded disk
$D\subset M$ whose boundary $\partial D$ is tangent to $kern~\alpha$
but $D$ itself is transversal to $kern~\alpha$ along $\partial D$. A
contact structure is said to be tight if it is not overtwisted. We
will call a contact form $\alpha$ overtwisted or tight if its
associated contact structure $kern~\alpha$ is overtwisted or tight
respectively. A codimension 1 foliation on a 3-manifold $M$ is said
to be taut if it is not the trivial foliation on $\text{\bf
S}^2\times\text{\bf S}^1$ and admits a riemannian metric $g$ for
which all of its leaves are minimal surfaces. Foliations defined by
closed, nonsingular 1-forms are taut. The following result can be
found on page 50 in \cite{ELT}.

\proclaim{Proposition 3.3} On a closed 3-dimensional manifold,
contact structures, $C^0$-close to a taut foliation, are tight.
\endproclaim
As a consequence of Proposition 3.3, one has the following:

\proclaim{Proposition 3.4} Suppose, on a closed 3-dimensional
manifold with a closed, nonsingular 1-form $\alpha_0$, that the
family of 1-forms $\alpha_t=\alpha_0+t\alpha$ is contact for all
$t>0$. Then $t\alpha$ and $\alpha_t$ are tight contact forms for all
$t>0$.
\endproclaim
\demo{Proof} For small $t>0$, $\alpha_t$ are tight contact forms by
Proposition 3.3. It follows from Proposition 3.2 that $\alpha_t$ and
$t\alpha$ are tight contact forms for all $t>0$. \quad\qed
\enddemo

We also prove the following theorem. \proclaim{Theorem 3.1} Let $M$
be a closed, 3-dimensional manifold with a closed, nonsingular
irrational 1-form $\alpha_0$. If $\alpha_t=\alpha_0+t\alpha$ is a
linear deformation of $\alpha_0$ into contact forms, then $\alpha$
is a quasi-regular, tight K-contact form.
\endproclaim

\demo{Proof} As in the proof of Proposition 3.2, for each $t$, let
$\beta_s=\alpha_t-s\alpha_0$ and define $X_s$ by $\beta_s(X_s)=0$
and $i_{X_s}d\beta_s=-\alpha_0$. Both $X_s$ and the Reeb vector
field $Z$ of $\alpha$ are tangent to the foliation defined by
$\alpha_0$ and whose leaves are all dense in $M$. It also follows
from the definitions that $X_s$ and $Z$ are pointwise linearly
independent. Let $L$ be a dense leaf of $\alpha_0$ with a periodic
$Z$ orbit on it. Such a periodic orbit exists by Proposition 1.1. By
Proposition 3.2, the isotopy $\Theta_s$ generated by $X_s$ preserves
the Reeb foliation defined by $Z$, therefore, all $Z$ orbits on $L$
are diffeomorphic, and hence periodic. Since $L$ is a dense leaf, it
follows that $M$ contains a dense subset of $Z$-periodic points.
Since the subset of all periodic points is closed in $M$,  we deduce
that all points in $M$ are periodic. From Proposition 1.2, we
conclude that $\alpha$ is a K-contact form. That $\alpha$ is tight
follows from Proposition 3.4. \quad\qed
\enddemo
The following corollary is an immediate consequence of Theorem 3.1
and Proposition 1.3.

\proclaim{Corollary 3.1} On the 3-dimensional torus $\text{\bf
T}^3$, no closed, nonsingular, irrational 1-form deforms linearly
into contact forms.
\endproclaim

\head 4. Classification of linearly deformable 1-forms\endhead
\noindent Motivated by Proposition 2.2, we observed in \cite{DA2}
that any nonsingular, closed, rational 1-form on a $\text{\bf T}^2$
bundle over $\text{\bf S}^1$ deforms linearly into contact forms. A
complete classification of these $\text{\bf T}^2$ bundle has been
presented in \cite{GEG}. Based on this classification, we complete
the result in \cite{DA2} as follows:

\proclaim{Theorem 4.1} Let $M$ be a closed $3$-dimensional manifold
which is a $\text{\bf T}^2$-bundle over $\text{\bf S}^1$. Then we
have the following:

i) If the monodromy of $M$ is not periodic, any closed non singular
$1$-form $\alpha_0$ on $M$ is  linearly deformable into contact
forms.

ii) If the monodromy of $M$ is periodic, a closed non singular
$1$-form $\alpha_0$ on $M$ is linearly deformable into contact forms
if and only if $Kern ~\alpha_0$ defines a fibration of $M$.
 \endproclaim

\demo{Proof}

 Let $M$ be a $\text{\bf T}^2$- bundle over the circle $\text{\bf S}^1$ and
$\alpha_{0}$ a closed nonsingular 1-form on $M$ with monodromy
matrix $A\in SL_{2}(\text{\bf Z})$. We consider $3$ different cases:

Case $1.$ The trace of $A$ satisfies $|tr(A)|\geq 3$: The manifold
$M$ is a left quotient of $Sol^{3}$, the solvable Lie group defined
as a split extension of $\text{\bf R}^2$ by $\text{\bf R}$. The
first Betti number of $M$ is $1$ and by Proposition 3.1, there exist
a diffeomorphism of $M$ isotopic to the identity of $M$ such that
$\varphi^{*}\alpha_{0}= a\pi^{*}d\theta$, where $a\in \text{\bf R}$,
$\pi$ is the $\text{\bf T}^2$-bundle map over $\text{\bf S}^1$, and
$d\theta$ the volume form of $\text{\bf S}^1$. Following
Dathe-Rukimbira (\cite{DA2}) $a\pi^{*}d\theta$ is linearly
deformable into contact forms and so is its diffeomorphic image
$\alpha_0$.

Case $2.$  The trace of $A$ satisfies $Tr(A)=2$, $M$ is then a left
quotient of $Nil^{3}$, the 3-dimensional Heisenberg  group. The
manifold $M$ is a non trivial $\text{\bf S}^{1}$-bundle over
$\text{\bf T}^2$ and its first Betti number is $2$. The left
invariant $1$-forms are generated by the $1$-forms $\omega_1$,
$\omega_2$, $\omega_3$ such: $d\omega_{1}= d\omega_{2}=0$,
$d\omega_{3}= -\omega_{1}\wedge \omega_{2}$. By Proposition 3.1, the
$1$-form $\alpha_0$ is isotopic to $a\omega_{1}+ b\omega_{2}$. The
family of $1$-forms $\alpha_{t}= a\omega_{1}+ b\omega_{2} +
t\omega_{3}$ is a linear deformation of $a\omega_{1}+ b\omega_{2}$
into contact forms. Since $\alpha_0$ is diffeomorphic with
$a\omega_1+b\omega_2$, it is also linearly deformable into contact
forms.

Case $3.$ The matrix $A$ is periodic: $M$ is a compact left
quotients of $\tilde{E_{2}}$, the universal cover of the group of
euclidian motions in the plane. There are 5 diffeomorphism classes
of such manifolds, only the 3-torus has first Betti number $>1$.
Given that on a closed 3-manifold with first Betti number equal to
1, any nonsingular, closed 1-form deforms linearly into contact
forms, case $3$ is reduced to showing that a closed, nonsingular
1-form on $\text{\bf T}^3$ deforms linearly into contact forms if
and only it is a fibration 1-form. From (\cite{DA2}), we know that
if $\alpha_0$ is rational, then it is linearly deformable into
contact forms. Let now $\alpha_0$ be an irrational, 1-form on
$\text{\bf T}^3$. Corollary 3.1 above says that $\alpha_0$ doesn't
deform linearly into contact forms, which completes the proof of
Theorem 4.1.\quad\qed
\enddemo

\refstyle{ABC}


\widestnumber\key{ABC}

\Refs

\ref\key{BLA} \by Blair, D.E. \book Riemannian geometry of contact
and symplectic manifolds, Progress in Mathematics 203 \publ
Birkhauser \publaddr Boston, Basel, Berlin \yr 2002
\endref

\ref\key{DA1} \by Dathe,H., Rukimbira, P. \paper Foliations and
contact structures \jour Advances in Geometry \vol 4, No. 1 \yr
2004\pages 75--81
\endref
\ref\key{DA2} \by Dathe, H., Rukimbira, P. \paper Fibration and
contact structures\jour Intern. Journ. of Math. and Math.
Sciences\vol 4 \yr 2005\pages 555--560
\endref

\ref\key{ELI} \by Eliashberg, Y.M. \paper Classification of
overtwisted contact structures\jour  Invent. Math.\vol 98\yr
1989\pages 623--637\endref

\ref\key{ELT} \by Eliashberg, Y.M , Thurston, W.P \book
Confoliations, University Lecture Series 13 \publ Amer. Math. Soc.
\yr 1998\endref \ref\key{LAB} \by Laudenbach, F., Blank, S. \paper
Sur l'isotopie des formes ferm\'ees en dimensions $3$\jour Invent.
Math. \vol 54\yr 1979\pages 103--177
\endref
\ref\key{GEG}\by Geiges, H., Gonzalo, J.\paper Contact circles on
3-manifolds\jour J. Diff. Geom. \vol 46\yr 1997\pages 236--286
\endref
\ref\key{RUK} \by Rukimbira, P.\paper Chern-Hamilton's Conjecture
and K-contactness\jour Houston Juor. Math.\vol 21, no. 4\yr
1995\pages 709--718
\endref
\ref\key{RU1} \by Rukimbira, P.\paper Some remarks on R-contact
flows \jour Ann. Glob. Anal. Geom.\vol 11\yr 1993\pages 165--171
\endref
\ref\key{TAU} \by Taubes, C.H. \paper The Seiberg-Witten equations
and the Weinstein conjecture\jour Geom. Topol.\vol  11\yr 2007\pages
2117--2202
\endref
\ref\key{TIS} \by Tischler, D. \paper On fibering certain foliated
manifolds over $\text{\bf S}^1$\jour Topology \vol 9\yr 1970\pages
153-154
\endref












\endRefs
\enddocument